\documentclass{amsart}
\usepackage{amssymb,latexsym,amscd}

\newtheorem{th:classif}[subsection]{Theorem}
\newtheorem{le:weakly}[subsection]{Lemma}
\newtheorem{le:cover}[subsection]{Lemma}
\newtheorem{le:homo}[subsection]{Lemma}
\newtheorem{le:theta}[subsection]{Lemma}
\newtheorem{pr:herm}[subsection]{Proposition}
\newtheorem{le:nottube}[subsection]{Lemma}
\newtheorem{le:z'}[subsection]{Lemma}
\newtheorem{le:adjoint}[subsection]{Lemma}
\newtheorem{pr:su(2n+1)}[subsection]{Proposition}
\newtheorem{le:su(2n+1)}[subsection]{Lemma}
\newtheorem{co:su(2n+1)}[subsection]{Corollary}
\newtheorem{le:action}[subsection]{Lemma}
\newtheorem{le:metric}[subsection]{Proposition}
\newtheorem{th:weakly}[subsection]{Theorem}
\newtheorem{th:newexamples}[subsection]{Theorem}

\theoremstyle{definition}
\newtheorem{de:weakly}[subsection]{Definition}
\newtheorem{de:homo}[subsection]{Definition}

\theoremstyle{remark}
\newtheorem{re:symmetric}[subsection]{Remark}

\begin{document}

\newcommand{\Proof}{\noindent {\em Proof.}\ }
\newcommand{\Qed}{~$\Box$}
\newcommand{\Real}{\mathbb{R}}
\newcommand{\Complex}{\mathbb{C}}
\newcommand{\Beginmatrix}{\left(\!\!\begin{array}}
\newcommand{\Endmatrix}{\end{array}\!\!\right)}
\newcommand{\nt}{\noindent}
\newcommand{\se}{\vspace{6pt}}
\newcommand{\fk}{\mathfrak{k}}
\newcommand{\fh}{\mathfrak{h}}
\newcommand{\fg}{\mathfrak{g}}
\newcommand{\fo}{\mathfrak{o}}
\newcommand{\fp}{\mathfrak{p}}
\newcommand{\fq}{\mathfrak{q}}
\newcommand{\fa}{\mathfrak{a}}
\newcommand{\fb}{\mathfrak{b}}
\newcommand{\fz}{\mathfrak{z}}
\newcommand{\fs}{\mathfrak{s}}
\newcommand{\fu}{\mathfrak{u}}
\newcommand{\di}{\displaystyle}
\newcommand{\tr}{\mathrm{Tr}}
\newcommand{\Dn}{D^1}
\newcommand{\Gn}{G^1}
\newcommand{\Kn}{K^1}
\newcommand{\An}{A_\nu}
\newcommand{\mn}{\mu_{\nu}}
\newcommand{\vol}{\mathrm{vol}}
\newcommand{\Vt}{\mathbb{V}_{\tau}}
\newcommand{\im}{\mathrm{Im}\, }
\newcommand{\ap}{\times_{\Phi}}

\title{Compact weakly symmetric spaces and spherical pairs}
\author{Hi\^{e}\'{u} Nguy\^{e}\~{n}}
\date{8/7/98}
\address{Department of Mathematics, 
Rowan University, Glassboro, NJ 08028, USA.}
\email{nguyen@rowan.edu}
\subjclass{53C35, 32M15}
\keywords{weakly symmetric spaces, spherical pairs, Gelfand pairs.} 

\begin{abstract}
Let $(G,H)$ be a spherical pair and assume that $G$ is a
connected compact simple Lie group and $H$ a closed subgroup
of $G$.  We prove in this paper
that the homogeneous manifold $G/H$ is weakly symmetric with
respect to $G$ and possibly an additional fixed isometry $\mu$.
It follows that M. Kr\"{a}mer's classification
list of such spherical pairs also becomes a classification list of
compact weakly symmetric spaces.
In fact, our proof involves a case-by-case study of all the
spherical pairs on Kr\"{a}mer's list.
\end{abstract} 

\maketitle
\setcounter{tocdepth}{1} \tableofcontents

\section{Introduction}

Let $M=G/H$ be a weakly symmetric space with respect to an
isometry group $G$ and possibly an additional fixed isometry
$\mu$.  A. Selberg was able to prove in 
\cite{se:weakly} that the space of $G$-invariant
differential operators on $M$ is commutative.  If $G$ is
connected, then this is equivalent to the property of the 
function space $L^1(H\!\setminus\! G/H)$ being commutative,
i.e. $(G,H)$ is a Gelfand pair (cf. \cite{th:diff} and
\cite{ng:gelfand}), or the property
that every unitary irreducible representations of $G$ contain
at most one $H$-fixed vector, i.e. $(G,H)$ is a spherical pair
and $H$ is a called a spherical subgroup of $G$ (cf. \cite{kr:spharische} 
and \cite{he:groups}, Ch. IV, Theorem 3.4).

It was proven by J. Lauret in \cite{la:comm} 
that the converse to the above
statement is false, namely that not all Gelfand pairs $(G,H)$ are
weakly symmetric.  Lauret's result involved constructing
generalized Heisenberg groups $N=G/H$ 
endowed with a modified Heisenberg-type metric and showing
that $N$ cannot be weakly symmetric 
with respect to the full isometry group $I(N)$.  However,
we prove in this paper that the converse does hold (in the
context of spherical pairs) if we assume $G$
to be a connected compact simple Lie group.  
Therefore, even in light
of Lauret's counterexamples, our result demonstrates
that for many of the compact spherical pairs $(G,H)$, 
proving weak symmetry of $M$ becomes a very easy task.

We now give a precise statement of our main result.

\begin{th:classif} \label{th:classif}
Let $(G,H)$ be a spherical pair with $G$ a connected 
compact simple Lie group and $H$ a closed subgroup of $G$.
Then the homogeneous manifold $M=G/H$ is 
weakly symmetric with respect to $G$
and possibly an additional fixed isometry $\mu$.
\end{th:classif}

Our strategy for proving this theorem is as follows: if $(G,H)$ is a spherical
pair, then Kr\"{a}mer in \cite{kr:spharische} has proven that $(G,H_0)$ is
also a spherical pair, where $H_0$ is the identity component
of $H$.  Furthermore, $(G,H_0)$ must appear on his classification list
given in the same paper. 
This allows us then to go through a case-by-case study of each pair appearing
on his list and prove that the corresponding homogeneous manifold $M_0=G/H_0$
is weakly symmetric with respect to $G$ and possibly an
additional fixed isometry $\mu$.  This is done by
understanding how the isotropy representation reverses tangent
vectors in each case.  
It follows that $M=G/H$ is also weakly symmetric with
respect to $G$ and $\mu$
as it has an even larger isotropy subgroup than $M_0$.

\se \nt {\it Acknowledgements.} The author wishes to 
thank W. Ziller for his generous help and guidance.  Many of the
arguments for the spherical pairs of constant
curvature and of Cayley-type 
actually resulted from joint work between us.  Furthermore,
he has recently informed the author that D.N. Akhiezer and E.B.
Vinberg have announced a similar result in
\cite{av:weakly} but within the context of spherical varieties.
In fact, they prove that Theorem \ref{th:classif} holds
more generally for real reductive algebraic groups $G$. 

\section{Preliminaries}

Let $M=G/H$ be a homogeneous Riemannian manifold, $G$ a transitive group
of isometries and $\mu$ a fixed 
isometry of
$M$ (not necessarily in $G$) satisfying $\mu G \mu^{-1} = G$ and $\mu^2 
\in G$.

\begin{de:weakly}
$M$ is {\it weakly symmetric} with respect to $G$ and $\mu$ if given any 
two points
$x$ and $y$ in $M$, there exists an element $g \in G$ such that $gx=\mu y$ 
and $gy = \mu x$.
\end{de:weakly}

We now present some results that will make it easier for us
to identify weakly symmetric spaces in terms of their
linear isotropy representations.  The first is a lemma that
characterizes weak symmetry as the reversal
of tangent vectors, first observed by W. Ziller
\cite{zi:weakly}.  

\begin{le:weakly} (cf. \cite{ng:weakly}, Lemma 2.2) \label{le:weakly}
Assume that $M$ is connected and $\mu$ fixes the point $z_o$.  Define 
$H$ to be the
isotropy subgroup of $G$ at $z_o$ and $T_{z_o}(M)$ to be the tangent 
space of $M$
at $z_o$.  Then $M$ is weakly symmetric with respect to $G$ and $\mu$ 
if and only
if given any tangent vector $v \in T_{z_o}(M)$,
there exists an element $h \in H$ such that $d(h\circ \mu)_{z_o}(v) = -v$. 
\end{le:weakly}

As weakly symmetric spaces are inherently homogeneous
manifolds, Lemma \ref{le:weakly} now allows us to also view weak symmetry of
$M$ in
terms of the adjoint representation of $H$.  We wish to develop
this relationship as a next step.
Let $G$ be a connected Lie group and $\theta$ an involutive
automorphism of $G$.  If $H$ be a compact $\theta$-invariant subgroup
of $G$, then there exists a reductive composition $\fg = \fh + \fq$
such that $\fq$ is $Ad_{G/H}(H)$-invariant.  Here, 
$Ad(H)=Ad_{G/H}(H)$ is the adjoint representation of $H$ on
$\fq\cong \fg/\fh$.

\begin{de:homo}
$(G,H,\theta)$ is called a {\it weakly symmetric triple} if 
given any element $X \in \fq$, there exists an element $h \in H$ such that
$(Ad(h)\circ d\theta)(X) = -X.$
\end{de:homo}

\begin{le:homo}
If $(G,H,\theta)$ is a weakly symmetric triple, then the
homogeneous manifold $M=G/H$, equipped with any $G$-invariant
Riemannian metric, is weakly symmetric
with respect to $G$ and the corresponding isometry $\mu$
induced from $\theta$.
\end{le:homo}

\begin{proof} Consider
the natural projection map $\pi: G \rightarrow G/H$.  
The automorphism $\theta$ induces an involutive diffeomorphism
$\mu$ of the homogeneous manifold $G/H$:
$$\mu(gH) = \theta(g)H, \ \ \ \ gH \in G/H.$$
Endow $M=G/H$ with any $G$-invariant Riemannian metric. 
Then $G$ becomes a transitive
group of isometries acting by left translations on right cosets of $H$. 
It follows that $\mu$ is an involutive isometry satisfying 
$\mu G \mu^{-1} = G$ and $\mu^{2} \in G$.  

Next, we identify $T_o(M)$, the tangent
space of $M = G/H$ at the origin $o=eK$, 
with $\fq$ via $\pi$.  It follows that the isotropy
action of $H$ on $T_o(M)$ is precisely the adjoint action of $H$ on $\fq$.
It becomes clear that $M$ is weakly symmetric with respect to $G$ and $\mu$ 
as a result of Lemma \ref{le:weakly}. 
\end{proof}

\begin{re:symmetric} It is trivial to see that if $G/K$ is a
symmetric space, where $K$ is an open subgroup of 
$G^{\theta}$ (the fixed
point set of the corresponding involution $\theta$), 
then $(G,K,\theta)$ is a weakly symmetric triple.  In fact,
the adjoint representation of $K$ is not necessary here to
reverse tangent vectors.
\end{re:symmetric}

\begin{le:cover} \label{le:cover} 
Let $\rho: \tilde{G}\rightarrow G$ be a covering
homomorphism between connected Lie groups, $\tilde{H} \subset G$ a connected
subgroup and $H=\rho(\tilde{H})$.  Then  
$\tilde{M} = \tilde{G}/\tilde{H}$ is
weakly symmetric with respect to $\tilde{G}$ if and only if
$M=G/H$ is weakly symmetric with respect to $G$.
\end{le:cover}

\begin{proof}  Both $\tilde{G}$ and $G$ have
isomorphic Lie algebras $\tilde{\fg}$ and $\fg$,
respectively, under $d\rho$.  Similarly, both $\tilde{H}$ 
and $H$ have isomorphic Lie algebras $\tilde{\fh}$ and
$\fh$, respectively.  It follows that $T_{\tilde{o}}(\tilde{M})
\cong \tilde{\fp}\cong \fp \cong T_o(M)$, where $\tilde{\fp}$
and $\fp$ come from the $Ad_{\tilde{G}/\tilde{H}}
(\tilde{H})$- and $Ad_{G/H}(H)$-invariant
decompositions $\fg = \fh + \fp$ and
$\tilde{\fg}=\tilde{\fh}+\tilde{\fp}$, respectively.  Furthermore, given
any $\tilde{h} \in \tilde{H}$ with $h=\rho(\tilde{h})$, the
following diagram commutes:
\[
\begin{CD}
\tilde{\fg} @>d\rho>> \fg \\
@VAd_{\tilde{G}/\tilde{H}}(\tilde{h})VV @VVAd_{G/H}(h)V \\
\tilde{\fg} @>>d\rho> \fg
\end{CD}
\]
Now, fix $\tilde{v} \in \tilde{\fp}$ and set $v=d\rho(\tilde{v})
\in \fp$.   Then $Ad_{\tilde{G}/\tilde{H}}
(\tilde{h})(\tilde{v})=-\tilde{v}$ for
some $\tilde{h}\in \tilde{H}$ if and only if $Ad_{G/H}(h)(v)=-v$
for some $h\in H$ with $h=\rho(\tilde{h})$.
\end{proof}

\section{Proof of Theorem}

As described earlier in the introduction,
we shall prove Theorem \ref{th:classif} by 
decomposing the spherical pairs on Kramer's classification list 
into six families (categorized below) and showing
weak symmetry for each family:

\se \nt \begin{tabular}{ll}
I. & Symmetric spaces, including $SO(8)/(SU(2)\cdot Sp(2))$.  \\
II. & $S^1$-bundles over hermitian symmetric spaces of
nontube type. \\ 
III. & $SU(2n+1)/Sp(n)$ and $SU(2n+1)/(Sp(n)\cdot U(1))$.  \\
IV. & Spaces of constant curvature: 
$G_2/A_2$, $SO(7)/G_2$ and $SO(8)/Spin(7)$. \\
V. & Spaces of Cayley-type: $SO(10)/(SO(2)\times Spin(7))$,
$SO(9)/Spin(7)$ and $SO(8)/G_2$. \\
VI. & $SO(2n+1)/U(n)$ and $Sp(n)/(Sp(n-1)\times U(1))$. \\
\end{tabular}

\subsection{I} If $G/K$ is a symmetric space, then it is
obviously weakly symmetric with respect to
the corresponding involutive isometry (the group $G$ is not needed
here). We note that the homogeneous manifold $SO(8)/(SU(2)\cdot Sp(2))$
appearing on Kramer's list is also a
symmetric space (cf. \cite{wz:symmetric}, Table 3, p. 325).

\subsection{II}
Let $G$ be a connected semisimple matrix Lie group with finite center and $K$
a maximal compact subgroup of $G$ such that $D=G/K$ is a hermitian
symmetric space.  Let $\fg = \fk + \fp$ be the Cartan decomposition of $\fg$
with respect to a Cartan involution $\sigma$, where $\fg$ and $\fk$ are
the Lie algebras of $G$ and $K$, respectively.  Then $\fk = \fk_s + \fz_\fk$,
where $\fk_s$ is the semisimple part of $\fk$ and $\fz_\fk$ the one-dimensional
center of $\fk$.  Let $K_s$ be the subgroup of $K$ with Lie algebra
$\fk_s$ and $Z_K$ the center of $K$.  Then $K_s$ is connected
and $K=K_sZ^0_K$, where $Z^0_K$ is the connected component of $Z_K$.

Fix a maximal abelian subspace $\fa$ of $\fp$ and let $A$ and $P$ be subgroups
of $G$ with Lie algebras $\fa$ and $\fp$, respectively.  According to 
Flensted-Jensen \cite{fj:spher1},
there exists an involutive automorphism $\theta$ satisfying:

\begin{le:theta} \label{le:theta}
(i) $\sigma\theta = \theta\sigma$,

\nt (ii) $\theta(a) = a^{-1}$ for all $a \in A$,

\nt (iii) $\theta(K_s) = K_s$ and $\theta(c) = c^{-1}$ for all $c \in Z^0_K$.

\se \nt At the Lie algebra level, these properties translate to

\se \nt (i$'$) $d\sigma d\theta = d\theta d\sigma$,

\nt (ii$'$) $d\theta(X) = -X$ for all $X \in \fa$,

\nt (iii$'$) $\theta(K_s) = K_s$ and $d\theta(X) = -X$ for all $X \in \fz_K$.
\end{le:theta}

\begin{pr:herm} \label{pr:herm} 
If $G/K$ is not
of tube type, then 
$(G,K_s,\theta)$ is a weakly symmetric triple.
\end{pr:herm}

In order to prove this proposition,
we shall need some finer structure theory about 
hermitian symmetric spaces due to Wolf and Kor\'{a}nyi \cite{kw:real}.
With their permission, we also 
present some unpublished results of theirs obtained in the early 
1980's \cite{wk:notes}.  

Extend $\fa$ to a Cartan subalgebra $\fh$ of $\fg$ and
let $Z^J \in \fz_\fk$ be the element which gives $\fp$ a complex structure $J$ 
corresponding to the root system for $(\fg,\fh)$ as described in \cite{kw:real}.
We then split $Z^J = Z^0 + Z'$, where $Z^0$ defines the
complex structure on the polydisc or the polysphere inside $G/K$ 
and $Z'$ is an element inside $\fk$ which centralizes $\fa$.

\begin{le:nottube} (\cite{kw:real})
$G/K$ is not of tube type if and only if $Z' \neq 0$.  
\end{le:nottube}

\begin{le:z'} (\cite{wk:notes})
Let $S' = \{\exp{tZ'}: t \in \Real\}$.  If $G/K$ is not of
tube type, then

\nt (i) $Z' \notin \fk_s$, 

\nt (ii) $K_s S' = K = S'K_s$. 
\end{le:z'}

\begin{proof}
(i) We first note that $Z^J \bot \fk_s$.  Then $Z^0\bot Z'$ because
$Z'$ centralizes $\fa$ and $Z^0 \in [\fa,J\fa] \subset ad(\fa)\fg$.
Were $Z' \in \fk_s$, then $0 = <Z^J,Z'> = <Z',Z'>$ and forces $Z' = 0$ 
(here $<\cdot, \cdot>$ is the Killing form on $\fg$).

\nt (ii) As $Z' \notin \fk_s$, the circle group $S'$ acts nontrivially on the
circle $K/K_s$.  Therefore, $S'$ is transitive there and hence $S'K_s = K$.
Use the map $k \mapsto k^{-1}$ to obtain $K_sS' = K$.
\end{proof}

\begin{le:adjoint} \label{le:adjoint}
If $G/K$ is not of tube type, then $Ad(K_s)(\fa) = \fp$. 
\end{le:adjoint}

\begin{proof}
It is a standard result that $Ad(K)(\fa)=\fp$.  Now write $K = K_s S'$
and use the fact that $S'$ centralizes $\fa$ to see the our lemma
immediately follows.
\end{proof}

\begin{proof} (second part of Prop. \ref{pr:herm}) 
Decompose $\fg = \fk_s
+\fz_\fk + \fp$.  Write any element $X \in \fz_\fk +\fp$ as $X = Z+ Y$ with
$Z\in \fz_\fk$ and $Y \in \fp$.  By Lemma \ref{le:adjoint}, we 
can write $Y = Ad(h)(W)$ with $h \in K_s$ and $W \in \fa$.
Then 
$$d\theta(Y) = Ad(\theta(h))(-W).$$
Set $k = h\theta(h)$.  Now use the fact that the adjoint action of $K$
commutes with $J$ and is trivial on $\fz_\fz$ to check
$$(Ad(k)\circ d\theta)(X) = -X.$$
This proves that $(G,K_s,\theta)$ is a weakly symmetric triple.
\end{proof}

\subsection{III}
We fix $G = SU(2n+1)$, $K = S(U(2n)\times U(1))$ and
$K_s = SU(2n)$
so that $G/K$ is complex
projective space $\mathbb{C}P^{2n}$.  
We choose a Cartan involution $\sigma$ so that $G$ and
$K$ have the following form as matrices (see Helgason \cite{he:diff}):
$$\begin{array}{l}
G = \{g \in SL(2n+1,\Complex): \ ^tgI_{p,q}\bar{g} = I_{p,q} \ \mathrm{and}
\ \det{g} = 1.\}, \ \ \ \ I_{p,q} = \Beginmatrix{cc} I_{2n} & 0 \\
0 & -1 \Endmatrix \\
K = \left\{ \Beginmatrix{cc} A & 0 \\ 0 & D \Endmatrix
\in G: A \in U(2n), D \in U(1) \ \mathrm{and} \ \det{A} \det{D} = 1. 
\right\} \\
K_s = \left\{ \Beginmatrix{cc} A & 0 \\ 0 & 1 \Endmatrix
\in G: A \in SU(2n).  \right\} \\
Z^0_K = \left\{ \Beginmatrix{cc} A & 0 \\ 0 & e^{-i2n\phi} \Endmatrix
\in K: A=\mathrm{Diag}(e^{i\phi},..., e^{i\phi}), \ \ \phi \in [0,2\pi).  \right\} \\
\end{array}$$

Let $\fa \subset \fp$ be the following:
$$
\fa = \left\{\Beginmatrix{cc} 0 & iw \\ -iw & 0 \Endmatrix 
: w \in \Real \right\}, 
\ \mathrm{where} \
\fp = \left\{\Beginmatrix{cc} 0 & W \\ \bar{W}^t & 0 \Endmatrix : W \ \mathrm{
complex}\ 2n \times 1 \ \mathrm{matrix}\right\}.$$
Then the corresponding involutive
automorphism $\theta$ of $G$ described in 
Lemma \ref{le:theta} becomes such that $\theta(g)$ is
complex conjugation of all the entries of $g$ as a matrix in $SL(n,\Complex)$.
As $G/K$ is not of tube type,
$(SU(2n+1), SU(2n), \theta)$ is a weakly symmetric triple
by Prop. \ref{pr:herm}. 

Let $H = Sp(n)$ be a maximal subgroup of $K_s = SU(2n)$ embedded
in $G$ as follows so that $K_s/H$ is a symmetric space:
$$H = \left\{ \Beginmatrix{cc} Sp(n) & 0 \\ 0 & 1 \Endmatrix
\right\}, \ \mathrm{where} \ Sp(n) = 
\left\{ \Beginmatrix{cc} A & B \\ -\bar{B} & \bar{A} \Endmatrix:
AA^* + BB* = I_n, \ AB^t = BA^t.\right\}$$

\begin{pr:su(2n+1)} \label{pr:su(2n+1)} $(SU(2n+1),Sp(n),\theta)$ is a 
weakly symmetric triple.
\end{pr:su(2n+1)}

\begin{proof} We first write $\fs\fu(2n+1) = \fs\fu(2n)+\fz_\fk + \fp$,
where $\fs\fu(2n) = \fk_s$.  Then
by writing $\fs\fu(2n) = \fs\fp(n) + \fq$, we can decompose it
further as $\fs\fu(2n+1) = \fs\fp(n) + \fq + \fz_\fk + \fp$.  
It now comes down to proving that given any element $X = V + Z + W \in
\fq + \fz_\fk + \fp$ with $V$, $Z$ and $W$ in $\fq$, $\fz_\fk$
and $\fp$, respectively, there exists an element $k \in Sp(n)$ such that
$(Ad(k)\circ d\theta)(X) = -X$.

The Lie subalgebras $\fs\fp(n)$ and $\fq$ of $\fs\fu(2n)$ are defined
as follows:
$$\begin{array}{rcl}
\fs\fp(n) & = & \left\{\Beginmatrix{ccc} V_1 & V_2 & 0 \\ -\bar{V}_2 & \bar{V}_1 
& 0 \\ 0 & 0 
& 0 \Endmatrix : V_1 \in \fu(n), V_2 \ n\times n \ \mathrm{complex \ matrix.} 
\right\} \\ \\
\fq & = & \left\{\Beginmatrix{ccc} V_1 & V_2 & 0 \\ \bar{V}_2 & -\bar{V}_1 
& 0 \\ 0 & 0 
& 0 \Endmatrix : V_1 \in \fs\fu(n), V_2 \in \fs\fo(n,\Complex). \right\}
\end{array}$$
We define a maximal abelian subalgebra $\fb$ of $\fq$ as follows:
$$\fb = \left\{\Beginmatrix{ccc} D & 0 & 0 \\ 0 & -\bar{D} & 0 \\ 0 & 0 
& 0 \Endmatrix  \in \fq: D=\mathrm{Diag}(id_1,..., id_n) \ \mathrm{with} \
d_1,...,d_n \in \Real \ \mathrm{and} \ \sum_{j=1}^n d_j = 0. \right\}$$

Notice that $(SU(2n),Sp(n),\theta|_{SU(2n)})$ 
is a weakly symmetric triple since
$\theta|_{SU(2n)}$, the restriction of $\theta$ to $SU(2n)$,
is also an involutive 
automorphism with $Sp(n)$ a $\theta$-invariant subgroup and
satisfies all the conditions of Lemma \ref{le:theta}.

\se As $SU(2n)/Sp(n)$ is a symmetric space,
we can write $V = Ad(h)(U)$, where $h \in Sp(n)$ and $U \in \fb$.  Then
as the adjoint action of $Sp(n)$ is trivial on $\fz_\fk$,
$$Ad(h^{-1})(X) = U + Z + W', \ \ \ \ W' \in \fp.$$
Now, let $N = N_{Sp(n)}(\fb)$ be the normalizer of $\fb$ in $Sp(n)$, i.e.
$$N(\fb) = \{ l \in Sp(n): Ad(l)(\fb) \subseteq \fb \}.$$
According to
Lemma \ref{le:su(2n+1)} below,
there exists an element $l \in N$ such that $Ad(l)(W') =
Y$ where $Y \in \fa$.  Hence,
$Ad(lh^{-1})(X) = U' + Z + Y$ with $U' \in \fb$.  In other words,
$$X = Ad(hl^{-1})(U' + Z + Y).$$
As $d\theta = -\mathrm{Id}$ on $\fb + \fz_\fk + \fa$, it follows that
$$d\theta(X) = Ad(\theta(hl^{-1}))(-U'-Z-Y).$$  
Setting $k = hl^{-1}\theta (hl^{-1}) \in Sp(n)$, we get
$$(Ad(k)\circ d\theta)(X) = -X.$$
\end{proof}
We also obtain as an obvious result the following corollary:

\begin{co:su(2n+1)} $(SU(2n+1),Sp(n)\cdot U(1),\theta)$ is also
a weakly symmetric triple.
\end{co:su(2n+1)}

\begin{le:su(2n+1)} \label{le:su(2n+1)}
$Ad_G(N)(\fa) = \fp$.
\end{le:su(2n+1)}

\begin{proof}
It is clear that $Ad_G(N)(\fa) \subset \fp$.  We introduce some
notation.  Let $E_{i,j}$ denote
the $n\times n$ matrix with 1 in the $ij$ entry and 0's everywhere
else, $F_{i,j}$ the $n\times n$ matrix obtained from
the identity matrix by switching the $i$-th row with
the $j$-th row, and $W_j$ the $(2n+1)\times 1$ matrix
with 1 in the $j$-th row and $0$'s everywhere else.  
Now, decompose $\fp = \sum_{j=1}^{2n} \fa_j$,
where $\fa_j$ is the complex subspace spanned by the matrix
$$\Beginmatrix{cc} 0 & W_j \\
W_j^t & 0 \Endmatrix.$$  
The following assertions can be easily verified:

\se \nt (i) Since $N$ contains the diagonal matrices in 
$Sp(n)$, it follows that $Ad(N)(\fa) \supset \fa_1$.

\se \nt (ii) Since $N$ contains matrices of the form
$$L_j = \Beginmatrix{ccc} F_{j,j+1} & 0 & 0 \\
0 & F_{j,j+1} & 0 \\
0 & 0 & 1 \Endmatrix, \ \ \ \ j =1,...,n,$$
it follows that $Ad(L_j)(\fa_j) = \fa_{j+1}$
and $Ad(L_j)(\fa_{n+j}) = \fa_{n+j+1}$.

\se \nt (iii) Since $N$ contains matrices of the form
$$\tilde{L}_j = \Beginmatrix{ccc} I_n - E_{j,j} & E_{j,j}
& 0 \\ -E_{j,j} & I_n - E_{j,j} & 0 \\
0 & 0 & 1 \Endmatrix, \ \ \ \ j=1,...,n,$$
it follows that $Ad(\tilde{L}_j)(\fa_j) = \fa_{n+j}$.

\se \nt From these assertions, we conclude that $Ad(N)(\fa)=\fp$.
This completes the proof of the lemma and Proposition \ref{pr:su(2n+1)}.
\end{proof}

\subsection{IV} Spaces of constant curvature:

\se \nt (i) $SO(8)/Spin(7)$: 
We have the double cover $Spin(8)/Spin(7)
\rightarrow SO(8)/Spin(7)$.  Since we have that $Spin(8)/Spin(7)\cong S^7$
is a space of constant curvature and 
the isotropy action is transitive on the unit tangent
sphere, it follows that $Spin(8)/Spin(7)$ is
weakly symmetric with respect $Spin(8)$. Hence, 
$SO(8)/Spin(7)$ is weakly symmetric with respect to
$SO(8)$ by Lemma \ref{le:cover}.

\se \nt (ii) $SO(7)/G_2$: Again, we have the double cover
$Spin(7)/G_2\rightarrow SO(7)/G_2$.  Now, $Spin(7)/G_2 \cong
S^7$ is a space of constant curvature and therefore weakly symmetric
with respect to $G_2$.  Hence, $SO(7)/G_2$ is
weakly symmetric with respect to $SO(7)$.

\se \nt (iii) $G_2/A_2$: From \cite{wo:spaces}, Prop. 8.12.7,
it is known that $G_2/A_2 = S^6$ is a space of constant
curvature and its isotropy action is transitive on the unit
tangent sphere. Hence, $G_2/A_2$ is weakly symmetric with respect to $G_2$.

\subsection{V} Spaces of Cayley-type:

\se \nt (i) $SO(10)/(SO(2)\times
Spin_{\pm}(7))$:  Following \cite{ta:go}, Theorem 5.4, we have the
decomposition
\[
\fs\fo(10)/(\mathfrak{spin}(7)\oplus
\fs\fo(2))= {\mathbb R}^7\oplus({\mathbb R}^8\otimes {\mathbb R}^2).
\]
Here, the adjoint representation of
$Spin_{\pm}(7)\times SO(2)$
on the factor $\mathbb{R}^7$ acts through the standard representation
of $Spin(7)$ on $\mathbb{R}^7$ (the
action of $SO(2)$ is trivial) and on the factor $\mathbb{R}^8\otimes
\mathbb{R}^2$ through the spin
representation of $Spin(7)$ on $\mathbb{R}^8$ and
by rotation of $SO(2)$ on $\mathbb{R}^2$.  Furthermore,
if $v \in\mathbb{R}^7$, then the isotropy
subgroup $(Spin(7)\times SO(2))_v = SU(4)\times SO(2)$, where
$SU(4)$ acts on 
$\mathbb{R}^8\cong \mathbb{C}^4$ by the standard action.

Therefore, given
any vector $v+w \in \mathbb{R}^7\oplus(\mathbb{R}^8\otimes
\mathbb{R}^2)$, there exists an element $k\in Spin(7)$ which
sends $v$ to its negative.  Next, we write
$w=\sum_{i=1}^pe_i\otimes f_j$, where
$e_i \in \mathbb{R}^8$ and $f_j \in \mathbb{R}^2$.
As $SO(2)$ acts by rotations, we use the rotation $e^{i\pi}$
to reverse each $f_j$ and therefore send $w$ to $-w$.
Together, the element $(k,e^{i\pi})$ in
$Spin(7)\times SO(2)$ will then send $v+w$ to $-(v+w)$. This
proves weak symmetry for $SO(10)/(SO(2)\times
Spin(7))$. 

\se \nt (ii) $SO(9)/Spin(7)$ and $SO(8)/G_2$:  
According to \cite{mu:exceptional}, there exists double
covers 
\[Spin(9)/Spin(7) \rightarrow SO(9)/Spin(7) \ \ \mathrm{and}
\ \ Spin(8)/G_2 \rightarrow SO(8)/G_2.\] 
But $Spin(9)/Spin(7)$ and
$Spin(8)/G_2$ are already known to be weakly symmetric with
respect to $Spin(9)$ and $Spin(8)$, respectively, from
\cite{zi:weakly}, pp. 357 and 361.

\subsection{VI} $SO(2n+1)/SU(n)$ and $Sp(n)/(U(1)\times
Sp(n-1))$: Again, these spaces have already been proven to
be weakly symmetric in \cite{zi:weakly}, pp. 359-360.

\se This concludes the proof of our main theorem.  
Table \ref{ta:weakly} is a summary of our results
with appropriate descriptions of weak
symmetry for each family.

\begin{table} \label{ta:weakly}
\begin{center}
\renewcommand{\arraystretch}{1.25}
\begin{tabular}{||c|ll|l|ll||} 
\multicolumn{6}{c}{Table \ref{ta:weakly} - Weakly symmetric
homogeneous manifolds of a connected compact simple Lie group} \\
\hline
\multicolumn{4}{||c|}{$M=G/H$ weakly symmetric} &
\multicolumn{2}{c||}{with respect to} \\ \hline
Family & $G$ \hspace{100pt} &  $H$ & for 
& & \\ \hline 
\renewcommand{\arraystretch}{1}
& & & & & \\ 
I & \multicolumn{2}{l|}{Symmetric spaces with involutive
isometry $\mu$} & & & $\mu$ \\
II & \multicolumn{2}{l|}{\parbox[t]{3.5in}{$S^1$-bundles over hermitian symmetric
spaces:}} & & & \\
& $SU(n+m)$ & $SU(m)\times SU(m)$ & $n>m\geq 1$ & $G$ \ \
and & $\mu$ \\
& $SO(2n)$ & $SU(n)$ & $n\geq 3$, $n$ odd & $G$ \ \ and &
$\mu$ \\
& $E_6$ & $D_5$ & & $G$ \ \ and & $\mu$ \\
III & $SU(2n+1)$ & $Sp(n)$ & $n\geq 1$ & $G$ \ \ and & $\mu$ \\
& $SU(2n+1)$ & $Sp(n)\cdot U(1)$ & $n\geq 1$ & $G$ \ \ and & $\mu$ \\ 
IV & \multicolumn{2}{l|}{Spaces of constant curvature:} & &
& \\
& $SO(8)$ & $Spin(7)$ & & $G$ & \\
& $SO(7)$ & $G_2$ & & $G$ & \\
& $G_2$ & $A_2$ & & $G$ & \\
V & \multicolumn{2}{l|}{Spaces of Cayley-type:} & & & \\  
& $SO(10)$ & $SO(2)\times Spin(7)$ & & $G$ & \\
& $SO(9)$ & $Spin(7)$ & & $G$ & \\
& $SO(8)$ & $G_2$ & & $G$ & \\
VI & $SO(2n+1)$ & $U(n)$ & $n \geq 2$ & $G$ & \\
& $Sp(n)$ & $Sp(n-1)\times U(1)$ & $n\geq 1$ & $G$ & \\
& & & & & \\ \hline 
\end{tabular}
\end{center}
\end{table}

\end{document}